\documentclass[a4paper,12pt]{amsart}
\usepackage{longtable,amsmath,amssymb,amsthm}
\title{On the Diophantine equation $x^2+q^{2m}=2y^p$}
\author{Sz. Tengely}
\address{Mathematical Institute\newline
 \indent University of Derecen\newline
 \indent P.O.Box 12\newline
 \indent 4010 Debrecen\newline
 \indent Hungary}
\email{tengely@math.klte.hu}
\keywords{Diophantine equations}
\subjclass[2000]{Primary 11D61; Secondary 11Y50}

\begin{document}
\newcommand{\Bound}{\mbox{Bound}}
\newcommand{\NE}{\mbox{NumofEq}}
\newcommand{\Root}{\mbox{Root}}
\newcommand{\pol}{\mbox{pol}}
\newcommand{\Red}{\mbox{Reduction}}
\newcommand{\ord}{\mbox{\rm ord}}
\newcommand{\lcm}{\mbox{\rm lcm}}
\newcommand{\sign}{\mbox{\rm sign}}
\newcommand{\ggd}{\mbox{\rm ggd}}

\newtheorem{thm}{Theorem}
\newtheorem{lem}{Lemma}
\newtheorem*{cor}{Corollary}
\newtheorem*{thm1}{Theorem}
\newtheorem*{lem1}{Lemma}
\newtheorem*{conj1}{Conjecture}
\theoremstyle{definition}
\newtheorem*{rem}{Remark}
\newtheorem*{acknowledgement}{Acknowledgement}
\bibliographystyle{plain}
\maketitle



\section{Introduction}
There are many results in the literature concerning the Diophantine equation
$$
Ax^2+p_1^{z_1}\cdots p_s^{z_s}=By^n,
$$
where $A,B$ are given non-zero integers, $p_1,\ldots,p_s$ are given primes and $n,x,y,$ $z_1,\ldots,z_s$ are integer unknowns with $n>2,x$ and $y$ coprime and non-negative, and $z_1,\ldots,z_s$ non-negative, see e.g. \cite{AM1}, \cite{AM2}, \cite{AM3}, \cite{AM4},
\cite{AM5}, \cite{AM6}, \cite{AM7}, \cite{B1}, \cite{C1}, \cite{C4}, \cite{Lu1}, \cite{Lu2}, \cite{Mig}, \cite{Mu1}, \cite{Mu2}, \cite{BP}.
Here the elegant result of Bilu, Hanrot and Voutier \cite{BHV} on the existence of primitive divisors of
Lucas and Lehmer numbers has turned out to be a very powerful tool.
Using this result Luca \cite{Lu2} solved completely the Diophantine equation $x^2+2^a3^b=y^n.$
Le \cite{Le} obtained necessary conditions for the solutions of the equation ${x^2 + p^2 = y^n}$ in positive integers ${x,y,n}$ with ${\gcd(x,y) = 1}$ and ${n>2}$. He also determined all solutions of this equation for ${p<100}$.
In \cite{BP} Pink considered the equation $x^2+(p_1^{z_1}\cdots p_s^{z_s})^2=2y^n,$ and gave an explicit
upper bound for $n$ depending only on $\max p_i$ and $s.$
The equation $x^2+1=2y^n$ was solved by Cohn \cite{C3}.
Pink and Tengely \cite{PT} considered the equation $x^2+a^2=2y^n.$ They gave an upper bound for the exponent $n$
depending only on $a,$ and completely resolved the equation with $1\leq a\leq 1000$ and $3\leq n\leq 80.$ In the present paper we study the equation $x^2+q^{2m}=2y^p$ where $m,p,q,x,y$ are integer unknowns with $m>0,$ $p$ and $q$ odd primes and $x$ and $y$ coprime. In Theorem \ref{c22_1} we show that all but finitely many solutions are of a special type. Theorem \ref{Baker} provides bounds for $p.$ Theorem \ref{c2fixy} deals with the case of fixed $y,$ we completely resolve the equation $x^2+q^{2m}=2\cdot 17^p.$ Theorem \ref{c2fixq} deals with the case of fixed $q.$ It is proved that if the Diophantine equation $x^2+3^{2m}=2y^p$ with $m>0$ and $p$ prime admits a coprime integer solution $(x,y),$ then $(x,y,m,p)\in\{(13,5,2,3),(79,5,1,5),(545,53,3,3)\}.$ It means that the equation $x^2+3^m=2y^p$ in copime integers is completely solved because solutions clearly do not exist when $m$ is odd.

\section{A finiteness result}
Consider the Diophantine equation
\begin{equation}\label{1}
x^2+q^{2m}=2y^p,
\end{equation}
where $x,y\in\mathbb{N}$ with $\gcd(x,y)=1, m\in\mathbb{N}$ and $p,q$ are odd primes and $\mathbb{N}$ denotes the set of positive integers.
Since the case $m=0$ was solved by Cohn \cite{C3} (he proved that the equation has only the solution $x=y=1$ in positive integers) we may assume without loss of generality that $m>0.$ If $q=2,$ then it follows from $m>0$ that $\gcd(x,y)>1,$ therefore we may further assume that $q$ is odd.
\begin{thm}\label{c22_1}
There are only finitely many solutions $(x,y,m,q,p)$ of \eqref{1} with $\gcd(x,y)=1, x,y\in\mathbb{N},$ such that $y$ is not a sum of two consecutive squares, $m\in\mathbb{N}$ and $p>3,q$ odd primes.
\end{thm}
\begin{rem}
The question of finiteness if $y$ is a sum of two consecutive squares is interesting.
The following examples show that very large solutions can exist.
\begin{center}
\tiny\begin{tabular}{|c|c|c|}
\hline
$y$ & $p$ & $q$\\
\hline
5 & 5 & 79\\
5 & 7 & 307\\
5 & 13 & 42641\\
5 & 29 & 1811852719\\
5 & 97 & 2299357537036323025594528471766399\\
13 & 7 & 11003\\
13 & 13 & 13394159\\
13 & 101 & 224803637342655330236336909331037067112119583602184017999\\
25 & 11 & 69049993\\
25 & 47 & 378293055860522027254001604922967\\
41 & 31 & 4010333845016060415260441\\
\hline
\end{tabular}
\end{center}
In these examples $m=1.$
\end{rem}
All solutions of \eqref{1} with small $q^m$ have been determined in \cite{T}.
\begin{lem}\label{small}
Let $q$ be an odd prime and $m\in\mathbb{N}\cup\{0\}$ such that $3\leq q^m\leq 501.$ If there exist $(x,y)\in\mathbb{N}^2$ with $\gcd(x,y)=1$ and an odd prime $p$ such that \eqref{1} holds, then
\begin{eqnarray*}
&&(x,y,q,m,p)\in\\
&&\bigl\{(3,5,79,1,5),(9,5,13,1,3),(13,5,3,2,3),(55,13,37,1,3),\\
&&(79,5,3,1,5),(99,17,5,1,3),(161,25,73,1,3),(249,5,307,1,7),\\
&&(351,41,11,2,3),(545,53,3,3,3),(649,61,181,1,3),(1665,113,337,1,3),\\
&&(2431,145,433,1,3),(5291,241,19,1,3),(275561,3361,71,1,3)\bigr\}.
\end{eqnarray*}
\end{lem}
\begin{proof}
This result follows from Corollary 1 in \cite{T}.
\end{proof}

We introduce some notation.
Put
\begin{equation}\label{delta4}
\delta_4=\begin{cases}
1 \mbox{ if } p\equiv 1\pmod{4},\\
-1 \mbox{ if } p\equiv 3\pmod{4}.
\end{cases}
\end{equation}
and
\begin{equation}\label{delta8}
\delta_8=\begin{cases}
1 \mbox{ if } p\equiv 1\mbox{ or }3\pmod{8},\\
-1 \mbox{ if } p\equiv 5\mbox{ or }7\pmod{8}.
\end{cases}
\end{equation}
Since $\mathbb{Z}[i]$ is a unique factorization domain, \eqref{1} implies the existence
of integers $u,v$ with $y=u^2+v^2$ such that
\begin{equation}\label{FG}
\begin{split}
x=\Re((1+i)(u+iv)^p)=:F_p(u,v),\\
q^m=\Im((1+i)(u+iv)^p)=:G_p(u,v).
\end{split}
\end{equation}
Here $F_p$ and $G_p$ are homogeneous polynomials in $\mathbb{Z}[X,Y].$
\begin{lem}\label{fac}
Let $F_p,G_p$ be the polynomials defined by \eqref{FG}. We have
\begin{eqnarray*}
(u-\delta_4 v)&|&F_p(u,v),\\
(u+\delta_4 v)&|&G_p(u,v).
\end{eqnarray*}
\end{lem}
\begin{proof}
This is Lemma 3 in \cite{T}.
\end{proof}
Lemma \ref{fac} and \eqref{FG} imply that there exists a $k\in\{0,1,\ldots,m\}$ such that
either
\begin{equation}\label{HP1}
\begin{split}
u+\delta_4 v=q^k,\\
H_p(u,v)=q^{m-k},
\end{split}
\end{equation}
or
\begin{equation}\label{HP2}
\begin{split}
u+\delta_4 v=-q^k,\\
H_p(u,v)=-q^{m-k},
\end{split}
\end{equation}
where $H_p(u,v)=\frac{G_p(u,v)}{u+\delta_4 v}.$

For all solutions with large $q^m$ we derive an upper bound for $p$ in case of $k=m$ in \eqref{HP1} or \eqref{HP2} and in case of $q=p.$
\begin{thm}\label{Baker}
If \eqref{1} admits a relatively prime solution $(x,y)\in\mathbb{N}^2$
then we have
\begin{eqnarray*}
&&p\leq 3803 \mbox{ if } u+\delta_4 v=\pm q^m, q^m \geq 503,\\
&&p\leq 3089 \mbox{ if } p=q,\\
&&p\leq 1309 \mbox{ if } u+\delta_4 v=\pm q^m, m\geq 40,\\
&&p\leq 1093 \mbox{ if } u+\delta_4 v=\pm q^m, m\geq 100,\\
&&p\leq 1009 \mbox{ if } u+\delta_4 v=\pm q^m, m\geq 250.\\
\end{eqnarray*}
\end{thm}
We shall use the following lemmas in the proof of Theorem \ref{Baker}.
The first result is due to Mignotte \cite[Theorem A.1.3]{BHV}.
Let $\alpha$ be an algebraic number, whose minimal polynomial over $\mathbb{Z}$ is
$A\prod_{i=1}^d(X-\alpha^{(i)}).$ The absolute logarithmic height of $\alpha$
is defined by
$$
h(\alpha)=\frac{1}{d}\left(\log|A|+\sum_{i=1}^d\log \max(1,|\alpha^{(i)}|)\right).
$$
\begin{lem}\label{Mig}
Let $\alpha$ be a complex algebraic number with $|\alpha|=1,$ but not a root of unity, and
$\log\alpha$ the principal value of the logarithm. Put $D=[\mathbb{Q}(\alpha):\mathbb{Q}]/2.$
Consider the linear form
$$
\Lambda=b_1 i\pi-b_2\log \alpha,
$$
where $b_1,b_2$ are positive integers. Let $\lambda$ be a real number satisfying $1.8\leq\lambda<4,$
and put
\begin{eqnarray*}
&&\rho=e^\lambda,\quad K=0.5\rho\pi+Dh(\alpha), \quad B=\max(13,b_1,b_2),\\
&&t=\frac{1}{6\pi\rho}-\frac{1}{48\pi\rho(1+2\pi\rho/3\lambda)},\quad T=\left(\frac{1/3+\sqrt{1/9+2\lambda t}}{\lambda}\right)^2,\\
&&H=\max\Bigl\{3\lambda,D\left(\log B+\log\left(\frac{1}{\pi\rho}+\frac{1}{2K}\right)-\log\sqrt{T}+0.886\right)+\\
&&+\frac{3\lambda}{2}+\frac{1}{T}\left(\frac{1}{6\rho\pi}+\frac{1}{3K}\right)+0.023\Bigr\}.
\end{eqnarray*}
Then
$$
\log|\Lambda|>-(8\pi T\rho\lambda^{-1}H^2+0.23)K-2H-2\log H+0.5\lambda+2\log\lambda-(D+2)\log 2.
$$
\end{lem}
The next result can be found as Corollary 3.12 at p. 41 of \cite{N}.
\begin{lem}\label{dep}
If $\Theta\in 2\pi\mathbb{Q},$ then the only rational values of the tangent and the cotangent
functions at $\Theta$ are $0,\pm 1.$
\end{lem}
\begin{proof}[Proof of Theorem \ref{Baker}]
Without loss of generality we assume that $p>1000$ and $q^m\geq 503.$
We give the proof of Theorem \ref{Baker} in the case $u+\delta_4 v=\pm q^m, q^m\geq 503,$ the proofs of the remaining four cases being analogous.
From $u+\delta_4 v=\pm q^m$ we get
$$
\frac{503}{2}\leq\frac{q^m}{2}\leq \frac{|u|+|v|}{2}\leq \sqrt{\frac{u^2+v^2}{2}}=\sqrt{\frac{y}{2}},
$$
which yields that $y\geq\frac{q^{2m}}{2}>126504.$
Hence
\begin{equation}\label{B}
\left|\frac{x+q^mi}{x-q^mi}-1\right|=\frac{2\cdot q^m}{\sqrt{x^2+q^{2m}}}\leq
\frac{2\sqrt{y}}{y^{p/2}}=\frac{2}{y^{\frac{p-1}{2}}}.
\end{equation}
We have
\begin{equation}\label{pth}
\frac{x+q^mi}{x-q^mi}=\frac{(1+i)(u+iv)^p}{(1-i)(u-iv)^p}=
i\left(\frac{u+iv}{u-iv}\right)^p.
\end{equation}
If $\left|i\left(\frac{u+iv}{u-iv}\right)^p-1\right|>\frac{1}{3}$ then
$6>y^{\frac{p-1}{2}}$, which yields a contradiction with $p>1000$ and $y>126504.$
Thus $\left|i\left(\frac{u+iv}{u-iv}\right)^p-1\right|\leq \frac{1}{3}.$
Since $|\log z|\leq 2|z-1|$ for $|z-1|\leq \frac{1}{3},$ we obtain
\begin{equation}\label{logineq}
\left|i\left(\frac{u+iv}{u-iv}\right)^p-1\right|\geq \frac{1}{2}\left|\log i\left(\frac{u+iv}{u-iv}\right)^p\right|.
\end{equation}

Suppose first that $\alpha:=\delta_4\left(\frac{u-iv}{-v+iu}\right)^{\sigma}$ is a root of unity for some $\sigma\in\{-1,1\}$. Then
$$
\left(\frac{u-iv}{-v+iu}\right)^{\sigma}=\frac{-2uv}{u^2+v^2}+\frac{\sigma(-u^2+v^2)}{u^2+v^2}i=\pm \alpha=\exp\left(\frac{2\pi ij}{n}\right),
$$
for some integers $j,n$ with $0\leq j\leq n-1.$ Therefore
$$
\tan\left(\frac{2\pi j}{n}\right)=\frac{\sigma(-u^2+v^2)}{-2uv}\in\mathbb{Q} \mbox{ or } (u,v)=(0,0).
$$
The latter case is excluded.
Hence, by Lemma \ref{dep}, $\frac{u^2-v^2}{2uv}\in\{0,1,-1\}.$
This implies that $|u|=|v|,$ but this is excluded by the requirement that the solutions $x,y$ of
\eqref{1} are relatively prime, but $y>126504.$
Therefore $\alpha$ is not a root of unity.

Note that $\alpha$ is irrational, $|\alpha|=1,$ and it is a root of the polynomial
$(u^2+v^2)X^2+4\delta_4 uvX+(u^2+v^2).$ Therefore $h(\alpha)=\frac{1}{2}\log y.$

Choose $l\in\mathbb{Z}$ such that $|p\log(i^{\delta_4}\frac{u+iv}{u-iv})+2l\pi i|$ is minimal, where logarithms have their principal values. Then $|2l|\leq p.$
Consider the linear form in two logarithms ($\pi i=\log(-1)$)
\begin{equation}\label{lambda}
\Lambda=2|l|\pi i-p\log \alpha.
\end{equation}

If $l=0$ then by Liouville's inequality and Lemma 1 of \cite{V},
\begin{equation}\label{Liou}
|\Lambda|\geq |p\log \alpha|\geq |\log \alpha|\geq 2^{-2}\exp(-2h(\alpha))\geq \\
\exp(-8(\log 6)^3h(\alpha)).
\end{equation}
From \eqref{B} and \eqref{Liou} we obtain
$$
\log 4-\frac{p-1}{2}\log y\geq \log|\Lambda|\geq -4(\log 6)^3\log y.
$$
Hence $p\leq 47.$ Thus we may assume without loss of generality that $l\neq 0.$

We apply Lemma \ref{Mig} with $\sigma=\sign(l), \alpha=\delta_4(\frac{u-iv}{-v+iu})^{\sigma}, b_1=2|l|$ and $b_2=p.$
Set $\lambda=1.8.$ We have $D=1$ and $B=p.$
By applying \eqref{B}-\eqref{lambda} and Lemma \ref{Mig} we obtain
\begin{equation*}\label{L}
\log 4-\frac{p-1}{2}\log y\geq \log|\Lambda|\geq -(13.16H^2+0.23)K-2H-2\log H-0.004.
\end{equation*}
We have
\begin{eqnarray*}
&&15.37677\leq K<9.5028+\frac{1}{2}\log y,\\
&&0.008633<t<0.008634,\\
&&0.155768<T<0.155769,\\
&&H<\log p+2.270616,\\
&&\log y>11.74803,
\end{eqnarray*}
From the above inequalities we conclude that
$p\leq 3803.$
\end{proof}

The following lemma gives a more precise description of the polynomial $H_p.$
\begin{lem}\label{polyH}
The polynomial $H_p(\pm q^k-\delta_4 v,v)$ has degree $p-1$ and
$$H_p(\pm q^k-\delta_4 v,v)=\pm\delta_82^{\frac{p-1}{2}}pv^{p-1}+q^kp\widehat{H}_p(v)+q^{k(p-1)},$$ where $\widehat{H}_p\in\mathbb{Z}[X]$ has degree $<p-1.$ The polynomial $H_p(X,1)\in\mathbb{Z}[X]$ is irreducible and
$$
H_p(X,1)=\prod_{\substack{k=0\\k\neq k_0}}^{p-1}\left(X-\tan\frac{(4k+3)\pi}{4p}\right),
$$
where $k_0=\left[\frac{p}{4}\right](p \mod 4).$
\end{lem}
\begin{proof}
By definition we have
\begin{equation}\label{defHp}
H_p(u,v)=\frac{G_p(u,v)}{u+\delta_4 v}=\frac{(1+i)(u+iv)^p-(1-i)(u-iv)^p}{2i(u+\delta_4 v)}.
\end{equation}
Hence
$$
H_p(\pm q^k-\delta_4 v,v)=\frac{(1+i)(\pm q^k+(i-\delta_4)v)^p-(1-i)(\pm q^k+(-i-\delta_4)v)^p}{\pm 2iq^k}.
$$
Therefore the coefficient of $v^p$ is $(1+i)(-\delta_4+i)^p+(1-i)(\delta_4+i)^p.$
If $\delta_4=1,$ then it equals
$-2(-1+i)^{p-1}+2(1+i)^{p-1}=-2(-4)^{\frac{p-1}{4}}+2(-4)^{\frac{p-1}{4}}=0,$ since $p\equiv 1\pmod{4}.$
If $\delta_4=-1,$ then it equals $(1+i)^{p+1}-(-1+i)^{p+1}=(-4)^{\frac{p+1}{4}}-(-4)^{\frac{p+1}{4}}=0.$
Similarly the coefficient of $v^{p-1}$ is $\pm\frac{(1+i)(\delta_4-i)^{p-1}-(1-i)(\delta_4+i)^{p-1}}{2i}p=\pm\delta_82^{\frac{p-1}{2}}p.$
It is easy to see that the constant is $q^{k(p-1)}.$ The coefficient of $v^{t}$ for $t=1,\ldots,p-2$ is $\pm\binom{p}{t}(q^k)^{p-t-1}c_t,$ where $c_t$ is a power of 2.
The irreducibility of $H_p(X,1)$ follows from the fact that $H_p(X-\delta_4,1)$
satisfies Eisenstein's irreducibility criterion. The last statement of the lemma is a direct consequence of Lemma 4 from \cite{T}.
\end{proof}
\begin{rem}
A deep conjecture, known as Schinzel's Hypothesis, says that 
if $P_1(X), ..., P_k(X)\in\mathbb{Z}[X]$ are irreducible polynomials such that no integer $l>1$ divides $P_i(x)$ for all integers $x$ for some $i\in\{1,\ldots,k\}$, then there exist infinitely many $x$ such that $P_1(x), ..., P_k(x)$ are simultaneously prime. Since $H_p(\pm 1-\delta_4 v,v)$ is irreducible having constant term $\pm 1,$ the Hypothesis implies that in case of $k=0,m=1$ there are infinitely many solutions of \eqref{1}.
\end{rem}
\begin{lem}\label{keq}
If there exists a $k\in\{0,1,\ldots,m\}$ such that
\eqref{HP1} or \eqref{HP2} has a solution $(u,v)\in\mathbb{Z}^2$ with $\gcd(u,v)=1$, then either $k=0$ or
$k=m, p\neq q$ or $(k=m-1,p=q).$
\end{lem}
\begin{proof}
Suppose $0<k<m.$ It follows from Lemma \ref{polyH} that $q\mid \pm\delta_82^{\frac{p-1}{2}}pv^{p-1}.$ If $q\neq p,$ we obtain that $q\mid v$ and $q\mid u,$ which is a contradiction with $\gcd(u,v)=1$. Thus $k=0$ or $k=m.$ If $p=q,$ then from Lemma \ref{polyH} and \eqref{HP1}, \eqref{HP2} we get
$$
\pm\delta_82^{\frac{p-1}{2}}v^{p-1}+p^k\widehat{H}_p(v)+p^{k(p-1)-1}=\pm p^{m-k-1}.
$$
Therefore $k=0$ or $k=m-1.$
\end{proof}

Now we are in the position to prove Theorem \ref{c22_1}.
\begin{proof}[Proof of Theorem \ref{c22_1}]
By Lemma \ref{keq} we have that $k=0,m-1$ or $k=m.$
If $k=0,$ then $u+\delta_4 v=\pm 1$ and $y$ is a sum of two consecutive squares.
If $k=m-1,$ then $p=q.$ Hence $u+\delta_4 v=\pm p^{m-1}$ which implies that $y\geq \frac{p^{2(m-1)}}{2}\geq \frac{p^2}{2}.$ From Theorem \ref{Baker} we obtain that $p\leq 3089.$
We recall that $H_p(u,v)$ is an irreducible polynomial of degree $p-1.$ Thus we have only finitely many Thue equations (if $p>3$)
$$
H_p(u,v)=\pm p.
$$
By a result of Thue \cite{Th} we know that
for each $p$ there are only finitely many integer solutions, which proves the statement.

Let $k=m.$ Here we have $u+\delta_4 v=\pm q^m$ and $H_p(\pm q^m-\delta_4 v,v)=\pm 1.$ If $q^m\leq 501$ then there are only finitely many solutions which are given in Lemma \ref{small}. We have computed an upper bound for $p$ in Lemma \ref{Baker} when $q^m\geq 503.$
This leads to finitely many Thue equations 
$$
H_p(u,v)=\pm 1.
$$
From Thue's result \cite{Th} follows that there are only finitely many integral solutions $(u,v)$ for any fixed $p,$
which implies the remaining part of the theorem.
\end{proof}
\section{Fixed $y$}
First we consider \eqref{1} with given $y$ which is not a sum of two consecutive squares. Since $y=u^2+v^2$ there are only finitely many possible pairs $(u,v)\in\mathbb{Z}^2.$ Among these pairs we have to select those for which
$u\pm v=\pm q^{m_0},$ for some prime $q$ and for some integer $m_0.$ Thus there are only finitely many pairs $(q,m_0).$ The method of \cite{T} makes it possible to compute (at least for moderate $q$ and $m_0$) all solutions of $x^2+q^{2m_0}=2y^p$ even without knowing $y.$
Let us consider the concrete example $y=17.$
\begin{thm}\label{c2fixy}
The only solution $(m,p,q,x)$ in positive integers $m,p,q,x$ with $p$ and $q$ odd primes of the equation $x^2+q^{2m}=2\cdot 17^p$ is $(1,3,5,99).$
\end{thm}
\begin{proof}
Note that 17 cannot be written as a sum of two consecutive squares.
From $y=u^2+v^2$ we obtain that $q$ is 3 or 5 and $m=1.$ This implies that 17 does not divide $x.$ We are left with the equations
\begin{eqnarray*}
&&x^2+3^2=2\cdot 17^p,\\
&&x^2+5^2=2\cdot 17^p.
\end{eqnarray*}
From Lemma \ref{small} we see that there is no solution with $q=3,m=1,y=17$ and the only solution in case of the second equation is $(x,y,q,m,p)=(99,17,5,1,3).$
\end{proof}
\section{Fixed $q$}
If $m$ is small, then one can apply the method of \cite{T} to obtain all solutions. Theorem \ref{Baker} provides an upper bound for $p$ in case $u+\delta_4 v=\pm q^m.$ Therefore it is sufficient to resolve the Thue equations
$$
H_p(u,v)=\pm 1
$$
for primes less than the bound. In practice this is a difficult job but in some special cases there exist methods which work, see \cite{BH}, \cite{BH1}, \cite{BHV}, \cite{H}.
Lemma \ref{subfield} shows that we have a cyclotomic field in the background just as in \cite{BHV}.
Probably the result of the following lemma is in the literature, but we have not 
found a reference. We thank Peter Stevenhagen for the short proof.
\begin{lem}\label{subfield}
For any positive integer $M$ denote by $\zeta_M$ a primitive $M$th root of unity.
If $\alpha$ is a root of $H_p(X,1)$ for some odd prime $p,$ then
$\mathbb{Q}(\zeta_p+\overline{\zeta}_p)\subset \mathbb{Q}(\alpha)\cong \mathbb{Q}(\zeta_{4p}+\overline{\zeta}_{4p}).$
\end{lem}
\begin{proof}
Since $\tan z=\frac{1}{i}\frac{\exp(iz)-\exp(-iz)}{\exp(iz)+\exp(-iz)},$ we can write
$\alpha=\tan(\frac{(4k+3)\pi}{4p})$ as
$$
\frac{1}{i}\frac{\zeta_{8p}^{4k+3}-\zeta_{8p}^{-4k-3}}{\zeta_{8p}^{4k+3}+\zeta_{8p}^{-4k-3}}=
-\zeta_{4}\frac{\zeta_{4p}^{4k+3}-1}{\zeta_{4p}^{4k+3}+1}\in\mathbb{Q}(\zeta_{4p}).
$$
Since it is invariant under complex conjugation, $\alpha$ is an element of $\mathbb{Q}(\zeta_{4p}+\overline{\zeta}_{4p}).$ We also know that $[\mathbb{Q}(\zeta_{4p}+\overline{\zeta}_{4p}):\mathbb{Q}]=[\mathbb{Q}(\alpha):\mathbb{Q}]=p-1,$ thus $\mathbb{Q}(\zeta_{4p}+\overline{\zeta}_{4p})\cong \mathbb{Q}(\alpha).$
The claimed inclusion follows from the fact that $\zeta_p+\overline{\zeta}_p$ can be expressed easily in terms of $\zeta_{4p}+\overline{\zeta}_{4p}.$
\end{proof}
It is important to remark that the Thue equations $H_p(u,v)=\pm 1$ do not depend on $q.$ 
By combining the methods of composite fields \cite{BH1} and non-fundamental units \cite{H} for Thue equations we may rule out some cases completely.
If the method applies it remains to consider the cases $u+\delta_4 v=\pm 1$ and $p=q.$ 
If $q$ is fixed one can follow a strategy to eliminate large primes $p.$
Here we use the fact that when considering the Thue equation
\begin{equation}\label{Thue}
H_p(u,v)=\pm 1.
\end{equation}
we are looking for integer solutions $(u,v)$ for which $u+\delta_4 v$ is a power of $q.$
Let $w$ be a positive integer relatively prime to $q,$ then the set $S(q,w)=\{q^m \mod w: m\in\mathbb{N}\}$ has $\ord_w(q)$ elements.
Let
\begin{eqnarray*}
&&L(p,q,w)=\\
&&\left\{s\in\{0,1,\ldots,\ord_w(q)\}: H_p(q^s-\delta_4 v,v)=1 \mbox{ has a solution modulo } w\right\}.
\end{eqnarray*}
We search for numbers $w_1,\ldots,w_N$ such that $\ord_{w_1}(q)=\ldots=\ord_{w_N}(q)=:w,$ say.
Then
$$
m_0\mod w\in L(p,q,w_1)\cap\ldots\cap L(p,q,w_N),
$$
where $m_0\mod w$ denotes the smallest non-negative integer congruent to $m$ modulo $w.$
Hopefully this will lead to some restrictions on $m.$
As we saw before the special case $p=q$ leads to a Thue equation
$H_p(u,v)=\pm p$ and the previously mentioned techniques may apply even for large primes. In case of $u+\delta_4 v=\pm 1$ one encounters a family of superelliptic equations $H_p(\pm 1-\delta_4 v,v)=\pm q^m.$ We will see that sometimes it is possible to solve these equations completely using congruence conditions only.
\vskip 5pt
From now on we consider \eqref{1} with $q=3,$ that is
\begin{equation}\label{13}
x^2+3^{2m}=2y^p.
\end{equation}

The equation $x^2+3=y^n$ was completely resolved by Cohn \cite{C2a}. Arif and Muriefah \cite{AM3} found all solutions of the equation $x^2+3^{2m+1}=y^n.$ There is one family of solutions, given by $(x,y,m,n)=(10\cdot 3^{3t},7\cdot 3^{2t}, 5+6t,3).$ Luca \cite{Lu1} proved that all solutions of the equation $x^2+3^{2m}=y^n$ are of the form $x=46\cdot 3^{3t}, y=13\cdot 3^{2t}, m=4+6t, n=3.$

\begin{rem}
We note that equation \eqref{13} with odd powers of 3 is easily solvable. From $x^2+3^{2m+1}=2y^p$ we get
$$
4\equiv 2y^p\pmod{8},
$$
hence $p=1$ which contradicts the assumption that $p$ is prime.
\end{rem}
Let us first treat the special case $p=q=3.$
By \eqref{FG} and Lemma \ref{fac} we have
\begin{eqnarray*}
x&=&F_3(u,v)=(u+v)(u^2-4uv+v^2),\\
3^m&=&G_3(u,v)=(u-v)(u^2+4uv+v^2).
\end{eqnarray*}
Therefore there exists an integer $k$ with $0\leq k\leq m,$ such that
\begin{eqnarray*}
u-v&=&\pm 3^k,\\
u^2+4uv+v^2&=&\pm 3^{m-k}.
\end{eqnarray*}
Hence we have
$$
6v^2\pm 6(3^k)v+3^{2k}=\pm 3^{m-k}.
$$
Both from $k=m$ and from $k=0$ it follows easily that $k=m=0.$ This yields the solutions $(x,y)=(\pm 1,1).$
If $k=m-1>0,$ then $3\mid 2v^2\pm 1.$ Thus one has to resolve the system of equations
\begin{eqnarray*}
u-v&=&-3^{m-1},\\
u^2+4uv+v^2&=&-3.
\end{eqnarray*}
The latter equation has infinitely many solutions parametrized by
\begin{eqnarray*}
&&u=\frac{-\varepsilon}{2}\left((2+\sqrt{3})^{t-1}+(2-\sqrt{3})^{t-1}\right),\\
&&v=\frac{\varepsilon}{2}\left((2+\sqrt{3})^{t}+(2-\sqrt{3})^{t}\right),
\end{eqnarray*}
where $t\in\mathbb{N},\varepsilon\in\{-1,1\}.$ Hence we get that
\begin{equation}\label{recur}
\frac{1}{2}\left((3+\sqrt{3})(2+\sqrt{3})^{t-1}+(3-\sqrt{3})(2-\sqrt{3})^{t-1}\right)=\pm 3^{m-1}.
\end{equation}
The left-hand side of \eqref{recur} is the explicit formula of the linear recursive sequence defined by $r_0=r_1=3, r_t=4r_{t-1}-r_{t-2}, t\geq 2.$ One can easily check that 
\begin{eqnarray*}
&&r_t\equiv 0\pmod{27}\Longleftrightarrow t\equiv 5\mbox{ or } 14\pmod{18},\\
&&r_t\equiv 0\pmod{17}\Longleftrightarrow t\equiv 5\mbox{ or } 14\pmod{18}.
\end{eqnarray*}
Thus $m=2$ or $m=3.$ If $m=2, k=1,$ then we obtain the solution $(x,y)=(13,5),$
if $m=3, k=2,$ then we get $(x,y)=(545,53).$ From now on we assume that $p>3.$


As we mentioned, sometimes it is possible to handle the case $k=0$ using congruence arguments only. In case of $q=3$ it works.
\begin{lem}\label{k0}
In case of $q=3$ there is no solution of \eqref{HP1} and \eqref{HP2} with $k=0.$
\end{lem}
\begin{proof}
We give a proof for \eqref{HP1} which also works for \eqref{HP2}.
In case of \eqref{HP1} if $k=0,$ then $u=1-\delta_4 v.$
Observe that by \eqref{defHp}
\begin{itemize}
\item if $v\equiv 0\pmod{3},$ then $H_p(1-\delta_4 v,v)\equiv
1\pmod{3},$
\item if $v\equiv 1\pmod{3}$ and $p\equiv 1\pmod{4},$ then
$H_p(1-\delta_4 v,v)\equiv
1\pmod{3},$
\item if $v\equiv 1\pmod{3}$ and $p\equiv 3\pmod{4},$ then
$H_p(1-\delta_4 v,v)\equiv
\pm 1\pmod{3},$
\item if $v\equiv 2\pmod{3}$ and $p\equiv 1\pmod{4},$ then
$H_p(1-\delta_4 v,v)\equiv
\pm 1\pmod{3},$
\item if $v\equiv 2\pmod{3}$ and $p\equiv 3\pmod{4},$ then
$H_p(1-\delta_4 v,v)\equiv
1\pmod{3}.$
\end{itemize}
Thus $H_p(1-\delta_4 v,v)\not\equiv 0\pmod{3}.$ Therefore there is no $v\in\mathbb{Z}$ such that $H_p(1-\delta_4 v,v)=3^m,$ as should be the case by \eqref{HP1} and \eqref{HP2}.
\end{proof}
Finally we investigate the remaining case, that is $u+\delta_4 v=3^m.$ We remark that $u+\delta_4 v=-3^m$ is not possible because from \eqref{HP2} and Lemma \ref{polyH} we obtain $-1\equiv H_p(-3^m-\delta_4 v,v)\equiv 3^{k(p-1)}\equiv 1\pmod{p}.$
\begin{lem}\label{cong}
If there is a coprime solution $(u,v)\in\mathbb{Z}^2$ of \eqref{HP1} with $q=3, k=m,$ then
$p\equiv 5\mbox{ or } 11\pmod{24}.$
\end{lem}
\begin{proof}
In case of $k=m$ we have, by \eqref{HP1} and Lemma \ref{polyH},
\begin{equation}\label{keqm}
H_p(3^m-\delta_4 v,v)=\delta_82^{\frac{p-1}{2}}pv^{p-1}+3^mp\widehat{H}_p(v)+3^{m(p-1)}=1.
\end{equation}
Therefore
$$
\delta_82^{\frac{p-1}{2}}p\equiv 1\pmod{3}
$$
and we get that $p\equiv 1,5,7,11\pmod{24}$.
Since by Lemma \ref{small} the only solution of the equation $x^2+3^{2m}=2y^p$ with $1\leq m\leq 5$
is given by $(x,y,m,p)\in\{(79,5,1,5),(545,53,3,3)\}$, we may assume without loss of generality that $m\geq 6.$
To get rid of the classes 1 and 7 we work modulo $243.$
If $p=8t+1,$ then from \eqref{keqm} we have
$$
2^{4t}(8t+1)v^{8t}\equiv 1\pmod{243}.
$$
It follows that $243|t$ and the first prime of the appropriate form is 3889 which is larger than the bound we have for $p$. If $p=8t+7,$ then
$$
-2^{4t+3}(8t+7)v^{8t+6}\equiv 1\pmod{243}.
$$
It follows that $t\equiv 60\pmod{243}$ and it turns out that $p=487$ is in this class, so we work modulo $3^6$ to show that the smallest possible prime is larger than the bound we have for $p$. Here we have to resolve the case $m=6$ using the method from \cite{T}. This value of $m$ is not too large so the method worked.
We did not get any new solution. Thus $p\equiv 5\mbox{ or } 11\pmod{24}.$
\end{proof}
\begin{thm}\label{kill}
There exists no coprime integer solution $(x,y)$ of $x^2+3^{2m}=2y^p$ with $m>0$ and $p<1000, p\equiv 5\pmod{24}$ or $p\in\{131, 251, 491, 971\}$ prime. 
\end{thm}
\begin{proof}
To prove the theorem we resolve the Thue equations \eqref{Thue}
for the given primes. In each case there is a small subfield, hence we can apply the method of \cite{BH1}. We wrote a PARI \cite{PARI} script to handle the computation. 
We note that if $p=659$ or $p=827,$ then there is a degree 7 subfield, but the regulator is too large to get unconditional result, the same holds for $p=419, 683, 947,$ in these cases there is a degree 11 subfield.
In the computation we followed the paper \cite{BH1}, but at the end we skipped the enumeration step. Instead we used the bound for $|x|$ given by the formula (34) at page 318.
The summary of the computation is in Table 1. 
\begin{table}[h]
\tiny
\begin{longtable}{|c|c|c|c|c|c|c|c|c|c|c|c|c|c|c|}
\caption[]{Summary of the computation (AMD64 Athlon 1.8GHz)}\\
\hline $p$      &	$X_3$   &	time & $p$      &	$X_3$   &	time &
$p$      &	$X_3$   &	time & $p$      &	$X_3$   &	time &
$p$      &	$X_3$   &	time \\
\endhead
\hline
\hline 29 & 4 &	1s  & 173 & 2 &	6s  & 317 & 2 &	13s  & 557 & 2 & 27s  & 797 & 2	& 45s\\

\hline 53 & 3 &	2s  & 197 & 2 &	7s  & 389 & 2 &	25s  & 653 & 2 & 33s  & 821 & 2	& 56s\\

\hline 101 & 2 & 3s & 251 & 2 &	14s & 461 & 2 & 22s  & 677 & 2 & 28s & 941 & 2 & 62s\\
 
\hline 131 & 2 & 6s & 269 & 2 &	14s & 491 & 2 & 25s  & 701 & 2 & 37s & 971 & 2 & 75s\\
 
\hline 149 & 2 & 7s & 293 & 2 & 10s  & 509 & 2 & 23s & 773 & 2 & 44s & & & \\
\hline
\end{longtable}
\end{table}
We obtained small bounds for $|u|$ in each case. It remains to find the integer solutions of the polynomial equations $H_p(u_0,v)=1$ for the given primes with $|u_0|\leq X_3.$ There is no solution for which $u+\delta v=3^m, m>0,$ and the statement follows.
\end{proof}

The Thue equations related to the remaining primes ($p<1000$) were solved by G. Hanrot.

\begin{thm}[G. Hanrot]\label{Hanrot}
There exists no coprime integer solution $(x,y)$ of $x^2+3^{2m}=2y^p$ with $m>0$ and $$p\in\{59, 83, 107, 179, 227, 347, 419, 443, 467, 563, 587, 659, 683, 827, 947\}.$$
\end{thm}
\begin{proof}
By combining the methods of composite fields \cite{BH1} and non-fundamental units \cite{H} all Thue equations were solved related to the given primes. The computations were done using PARI. Most of the computation time is the time for $p-1$
LLL-reductions in dimension 3 on a lattice with integer entries of
size about the square of the Baker bound. The numerical precision required in the worst case ($p=587$) is 7700. The summary of the computation is in Table 2. 
We got small bounds for $|u|$ in each case. 
There is no solution for which $u+\delta v=3^m, m>0,$ and the statement follows.
\begin{table}[h]
\tiny
\begin{longtable}{|c|c|c|c|c|c|c|c|c|}
\caption[]{Summary of the computation (AMD Opteron 2.6GHz)}\\
\hline $p$      &	$X_3$   &	time & $p$      &	$X_3$   &	time &
$p$      &	$X_3$   &	time\\
\endhead
\hline
\hline 59	&	47	&	2s & 347	&	186	&	33m &
587	&	279	&	248m\\ 
\hline 83	&	62	&	9s & 419	&	216	&	67m &
659	&	1	&	3s\\
\hline 107	&	74	&	23s & 443	&	2	&	5s &
683	&	2	&	7s\\
\hline 179	&	111	&	2m29s & 467	&	233	&	102m &
827	&	2	&	4s\\
\hline 227	&	134	&	6m13s & 563	&	270	&	211m &
947	&	2	&	10s\\
\hline
\end{longtable}
\end{table}
\end{proof}


We recall that Cohn \cite{C3} showed that the only positive integer solution of $x^2+1=2y^p$ is given by $x=y=1.$

\begin{thm}\label{c2fixq}
If the Diophantine equation $x^2+3^{2m}=2y^p$ with $m>0$ and $p$ prime admits a coprime integer solution $(x,y),$ then 
$(x,y,m,p)=(13,5,2,3), (79,5,1,5), (545,53,3,3).$
\end{thm}
\begin{proof}
We will provide lower bounds for $m$ which contradict the bound for $p$ provided by Theorem \ref{Baker}.
By Theorem \ref{Baker} we have $p\leq 3803$ and by Lemma \ref{cong} we have $p\equiv 5\mbox{ or } 11\pmod{24}.$
We compute the following sets for each prime $p$ with $1000\leq p\leq 3803, p\equiv 5\mbox{ or } 11\pmod{24}:$
\begin{eqnarray*}
&&A5=L(p,3,242),\\
&&A16=L(p,3,136)\cap L(p,3,193)\cap L(p,3,320)\cap L(p,3,697),\\
&&A22=L(p,3,92)\cap L(p,3,134)\cap L(p,3,661),\\
&&A27=L(p,3,866)\cap L(p,3,1417),\\
&&A34=L(p,3,103)\cap L(p,3,307)\cap L(p,3,1021),\\
&&A39=L(p,3,169)\cap L(p,3,313),\\
&&A69=L(p,3,554)\cap L(p,3,611).
\end{eqnarray*}
In case of $A5$ we have $\ord_{242}3=5,$ hence this set contains those congruence classes modulo 5 for which \eqref{13} is solvable, similarly in case of the other sets.
How can we use this information?
Suppose it turns out that for a prime $A5=\{0\}$ and $A16=\{0\}.$ Then we know that $m\equiv 0\pmod{5\cdot 16}$ and Theorem \ref{Baker} implies $p\leq 1309.$ If the prime is larger than this bound, then we have a contradiction. In Table 3 we included those primes for which we obtained a contradiction in this way.
\begin{table}[h]
\tiny
\begin{longtable}{|c|c|c|c|c|c|c|c|c|c|}
\caption[]{Excluding some primes using congruences.}\\
\hline
$p$ & $\mod$ & $p$ & $\mod$ & $p$ & $\mod$ & $p$ & $\mod$ & $p$ & $\mod$\\
\endhead
\hline
1013&	16,27	&	1571&	5,22	&	1973&	16,22	&	2357&	16,22	&	3011&	5,22\\
1109&	16,22	&	1613&	16,22	&	1979&	16,22	&	2459&	16,22	&	3203&	16,22\\
1181&	16,22	&	1619&	16,22	&	2003&	16,22	&	2477&	16,22	&	3221&	16,22\\
1187&	16,22	&	1667&	16,22	&	2027&	16,22	&	2531&	5,22	&	3323&	16,22\\
1229&	16,22	&	1709&	16,22	&	2069&	16,22	&	2579&	16,22	&	3347&	16,22\\
1259&	16,22	&	1733&	16,22	&	2099&	16,22	&	2693&	16,22	&	3371&	5,22\\
1277&	16,22	&	1787&	16,22	&	2141&	16,22	&	2741&	16,27	&	3413&	16,22\\
1283&	16,22	&	1811&	5,22	&	2237&	16,22	&	2861&	16,22	&	3533&	16,22\\
1307&	16,22	&	1877&	16,27	&	2243&	16,22	&	2909&	16,22	&	3677&	16,22\\
1493&	16,22	&	1931&	5,22	&	2309&	16,27	&	2957&	16,22	&	3701&	16,22\\
1523&	16,22	&	1949&	16,22	&	2333&	16,22	&	2963&	16,22	&	&\\
\hline
\end{longtable}
\end{table}
In the columns $\mod$ the numbers $n$ are stated for which sets $An$ were used for the given prime. It turned out that only 4 sets were needed. In case of $5,22$ we have $m\geq 110, p\leq 1093,$ in case of $16,22$ we have $m\geq 176, p\leq 1093$ and in the case $16,27$ we have $m\geq 432, p\leq 1009.$
\begin{table}[h]
\tiny
\begin{longtable}{|c|c|c|c|c|c|c|c|c|}
\caption[]{Excluding some primes using CRT.}\\
\hline
$p$ 	&	$r_m$	&	$CRT$	&	$p$ &	$r_m$	&	$CRT$	& 	$p$  & $r_m$  &	$CRT$\\
\endhead
\hline
1019	&	384	&	5,16,27 &	2267 &	448	&	5,16,69	& 	3389 &	170 & 5,27,34\\
1061	&	176	&	5,16,39	&	2339 &	208	&	5,16,39	& 	3461 &	116 &	5,16,39\\
1091	&	580	&	5,16,27 &	2381 &	44	&	5,27,34	& 	3467 &	336 &	5,16,27\\
1163	&	586	&	5,27,34	&	2411 &	180	&	5,16,27 & 	3491 &  850 &	5,27,34\\
1301	&	416	&	5,16,39	&	2549 &	320	&	5,16,27 & 	3539 &	112 &	5,16,39\\
1427	&	270	&	5,27,34	&	2699 &	640	&	5,16,69	& 	3557 &	176 &	5,16,39\\
1451	&	340	&	5,16,27 &	2789 &	204	&	5,27,34	& 	3581 &	150 &	5,27,34\\
1499	&	112	&	5,16,39	&	2819 &	352	&	5,16,27 & 	3659 &	112 &	5,16,39\\
1637	&	121	&	5,27,34	&	2837 &	131	&	5,27,34	& 	3779 &	72  &	5,27,34\\
1901	&	304	&	5,16,39	&	2843 &	136	&	5,27,34	& 	3797 &	416 &	5,16,39\\
1907	&	102	&	5,27,34	&	3083 &	340	&	5,27,34	& 	3803 &	136 &	5,27,34\\
1997	&	170	&	5,27,34	&	3251 &	580	&	5,16,27 & 	     &	    &	\\
2213	&	170	&	5,27,34	&	3299 &	64	&	5,16,39	& 	     &	    &	\\
\hline
\end{longtable}
\end{table}
We could not exclude all primes using the previous argument, but there is an other way to use the computed sets. We can combine the available information by means of the Chinese remainder theorem.
Let $CRT([a5,a16,a39],[5,16,39])$ be the smallest non-negative solution of the system of congruences
\begin{eqnarray*}
&&m\equiv a5\pmod{5}\\
&&m\equiv a16\pmod{16}\\
&&m\equiv a39\pmod{39},
\end{eqnarray*}
where $a5\in A5, a16\in A16$ and $a39\in A39.$ Let $r_m$ be the smallest non-zero element of the set
$
\{CRT([a5,a16,a39],[5,16,39]): a5\in A5, a16\in A16, a39\in A39\},
$
In Table 4 we included the values of $r_m$ and the numbers related to the sets $A5-A69.$
We see that $m\geq r_m$ in all cases.
For example, if $p=1019$ then $m\geq 384,$ and Theorem \ref{Baker} implies $p\leq 1009,$ which is a contradiction.
For $p=2381$ we used $A5,A27$ and $A34,$ given by $A5=\{0,1,4\},A27=\{0,14,15,17\},A34=\{0,10\}.$
Hence
\begin{eqnarray*}
&&\{CRT([a5,a27,a34],[5,27,34]): a5\in A5, a27\in A27, a34\in A34\}=\\
&&=\{ 0, 44, 204, 476, 486, 554, 690, 986, 1394, 1404, 1836, 1880, 1904,\\
&& 2040, 2390,2526, 2754, 3230, 3240, 3444, 3716, 3740, 3876, 4226 \}.
\end{eqnarray*}
The smallest non-zero element is 44 (which comes from $[a5,a27,a34]=[4,17,10]$), therefore $m\geq 44$ and $p\leq 1309,$ a contradiction. In this way all remaining primes $>1000$ can be handled.
We are left with the primes $p<1000, p\equiv 5 \mbox{ or } 11 \pmod{24}$ 
They are mentioned in Theorem \ref{kill} and in Theorem \ref{Hanrot}.
\end{proof}

\begin{acknowledgement}
The author wish to thank Guillaume Hanrot for the computations related to Theorem \ref{Hanrot} and for giving some hints how to modify the PARI code, which was used in \cite{BHV}, to make the computations necessary for Theorem \ref{kill}.
Furthermore, we would like to thank Robert Tijdeman for his valuable remarks and suggestions, Peter Stevenhagen for the useful discussions on algebraic number theory, and for the proof of Lemma \ref{subfield} and K\'alm\'an Gy\H{o}ry for calling our
attention to Schinzel's Hypothesis.

\end{acknowledgement}

\bibliography{all}

\end{document}